# An aperiodic tiling made of one tile, a triangle.

Vincent Van Dongen, PhD

November 6, 2021

vincent.vandongen@gmail.com

## Abstract

How many different tiles are needed at the minimum to create aperiodicity? Several tilings made of two tiles were discovered, the first one being by Penrose in the 1970's. Since then, scientists discovered other aperiodic tilings made of two tiles, including the square-triangle one, a tiling that has been particularly useful for the study of dodecagonal quasicrystals and soft matters. An open problem still exists: Can one tile be sufficient to create aperiodicity? This is known as the 'ein-stein' problem. We present in this paper an aperiodic tiling made of one single tile: an isosceles right triangle. The tile itself is not aperiodic and therefore not a solution to the ein-stein problem but we present a set of substitution rules on the same tile that forces the tiling to be aperiodic. This paper presents its construction rules that proves its aperiodicity. We also show that this tiling offers an underlying dodecagonal structure close to the one of square-triangle tiling.

## About this paper

Can one tile fully cover the plane in an aperiodic manner and only aperiodically? This problem is known as the *'ein-stein'* problem. Note that 'ein stein' means 'one tile' in German.

The tiling presented here is made of a single tile that is an isosceles right triangle. This is not an aperiodic tile as it can tile periodically but our proposed tiling made with this tile is aperiodic. This document presents in detail non-overlapping substitution rules for generating the aperiodic tiling. The figure below shows the tiling at different stages of its construction. Note that we intentionally left the triangles *'blank'*, i.e. without any color nor pattern on them. This illustrates the fact that the aperiodicity of this tiling is due to the positioning of the tiles and nothing else. In our case, the triangle can only take four orientations and their alignment is edge-to-edge.

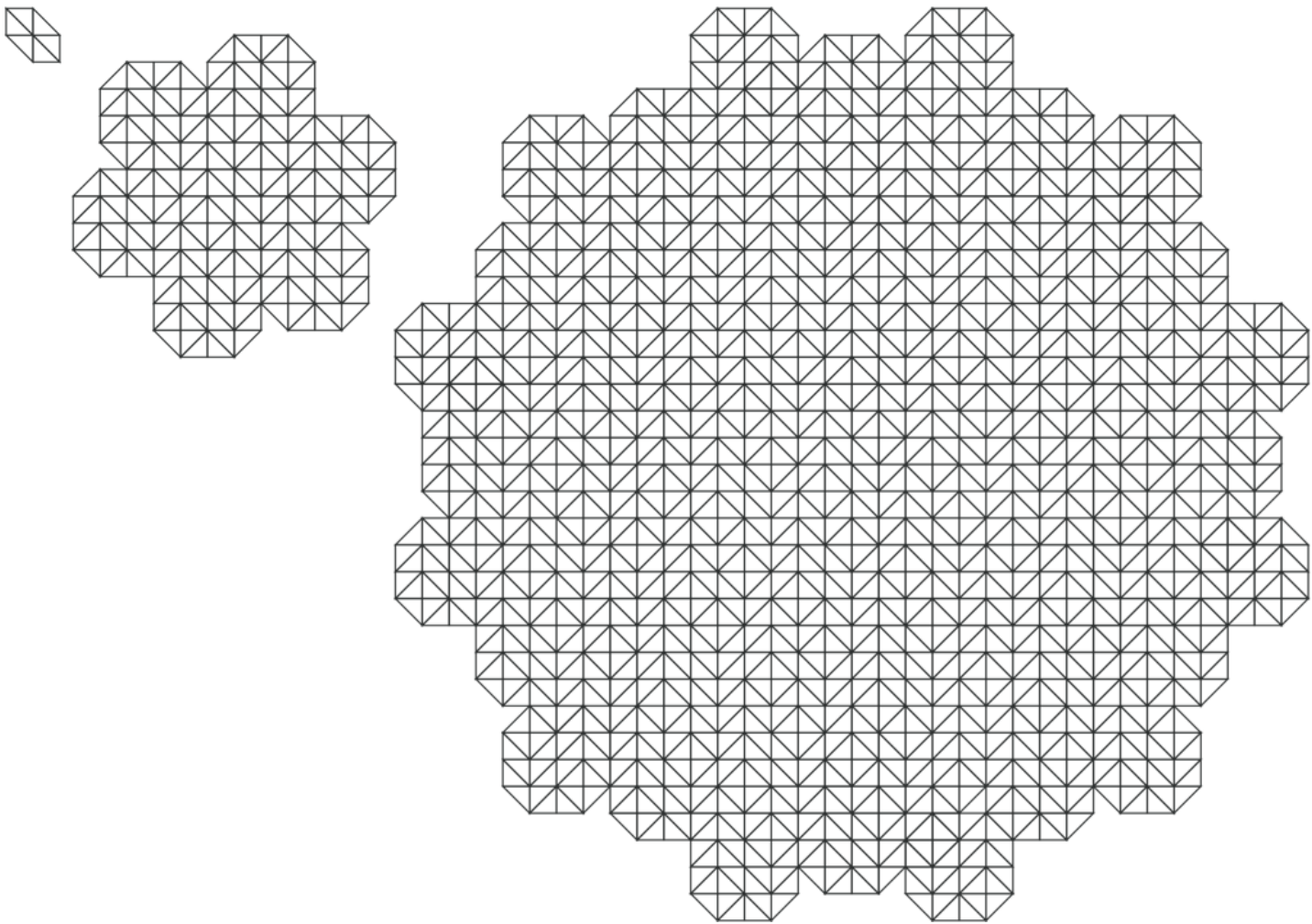

*Figure 1: An aperiodic tiling made of one tile, an isosceles right triangle.*

The tiling is aperiodic in the sense that its pattern keeps changing in a structured manner. The presentation of this new tiling is organized as follows. In the next section, we recall the concepts of aperiodic tiling and some well-known examples. In particular, the square-triangle tiling (M. Baake, 1992) is presented. In the next section, we show how to create a one-triangle aperiodic tiling. We then show its underlying dodecagonal structure, and its link with the square-triangle aperiodic tiling.

## Introduction on aperiodic tiling

These definitions are taken from (Wikipedia, Tesselation, s.d.):

A **tessellation** or tiling of a flat surface is the covering of a plane using one or more geometric shapes, called tiles, with no overlaps or no gaps.

An **edge-to-edge** tiling is any polygonal tessellation where adjacent tiles only share one full side, i.e., no tile shares a partial side or more than one side with any other tile. In an edge-to-edge tiling, the sides of the polygons and the edges of the tiles are the same.

A tiling is **periodic** when the tiling is made of the same pattern that repeats. Else, it is non-periodic.

The simplest periodic tiling to think of is the square grid. Squares are simply aligned edge-to-edge to pave the plane. In this case, the pattern that repeats is simply the square along its two edges and the period of the tile is the size of the edge.

Other simple tilings can be created from that grid. Indeed, each square of the grid can be divided in two along one of its two diagonals. If all diagonals have the same orientation, then the tiling is clearly periodic. But if the diagonals change orientation randomly, then the tiling is non-periodic. Here below is a comparison of three tilings, all made with the same tile: an isosceles right triangle. The one on the left is periodic as the same pattern repeats. The one in the middle is non periodic as no pattern can really be observed. The one on the right is the aperiodic tiling that will be presented in details in this paper.

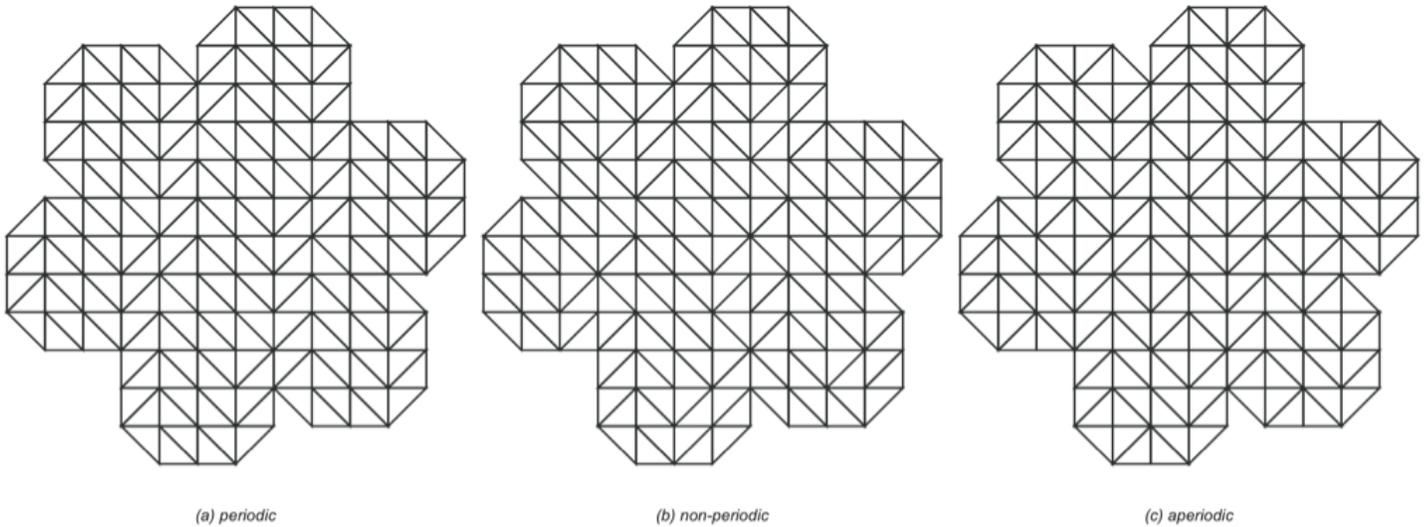

(a) periodic  (b) non-periodic  (c) aperiodic

Some non-periodic tilings can be structured and not aperiodic. Examples of these include the so-called Rep-tiles. These are self-replicating tessellations. An example is the so-called Sphinx tiling. See the figure below, taken from Wikipedia (Sphinx tiling).

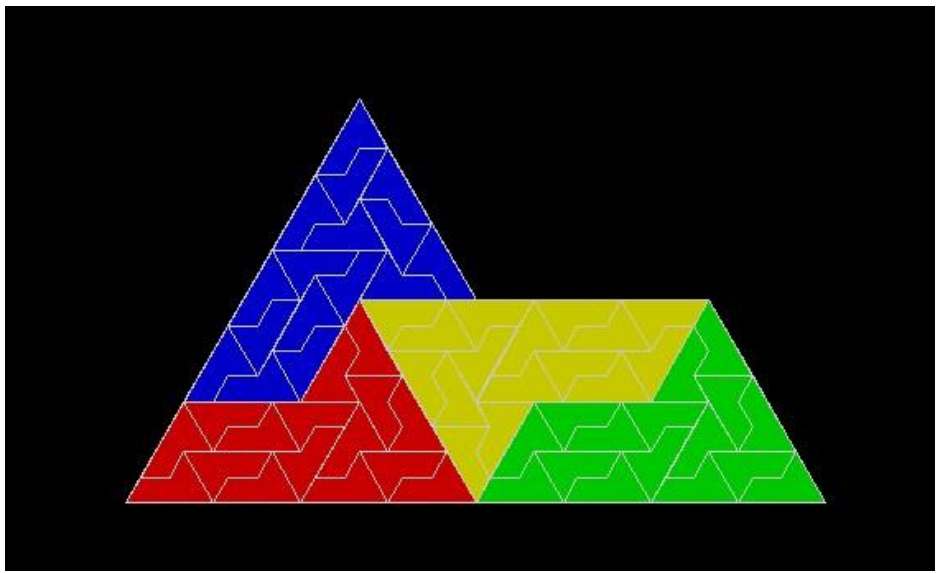

*Figure 2: The Sphinx tiling is an example of a self-replicating tessellation called Rep-Tile. Made of a single-tile, it generates a large number of different patterns but not an infinite number. Hence, this is not considered to be a one-tile aperiodic tiling. But it is non-periodic.*

Although the Sphinx Tiling generates many different patterns, this number is not infinite. Hence it is defined to be limitperiodic (Sphinx).

**Aperiodic** tiling is a special case of non-periodic tiling. Here is a definition taken from MathWorld (Weisstein): *An aperiodic tiling is a non-periodic tiling in which arbitrarily large periodic patches do not occur. A set of tiles is said to be aperiodic if they can form only non-periodic tilings. The most widely known examples of aperiodic tilings are those formed by Penrose tiles.*

Aperiodic tiling became popular when their link with quasi-periodic materials was established in the eighties. At this time, Professor Dan Shechtman, who received the Nobel Prize of Chemistry in 2011, had seen the existence of quasi-periodic materials through TEM, transmission electron microscopy (Shechtman). A link was made with the first 2-tile aperiodic tiling discovered by Penrose that exhibits a fivefold symmetry (Penrose, 1979).

Aperiodic tilings exhibit structures in the sense that they contain patterns. Each pattern is unique in a zone but the pattern can be found in a structured manner an infinite number of times in the tiling. But no single pattern of an aperiodic tiling ever covers the entire plane. In fact, an infinite number of patterns exist in the aperiodic tiling that we are interested in and each of them only partially covers the plane.

Let us recall the square-triangle tiling (Schlottmann). It consists of two tiles: an equilateral triangle and a square.

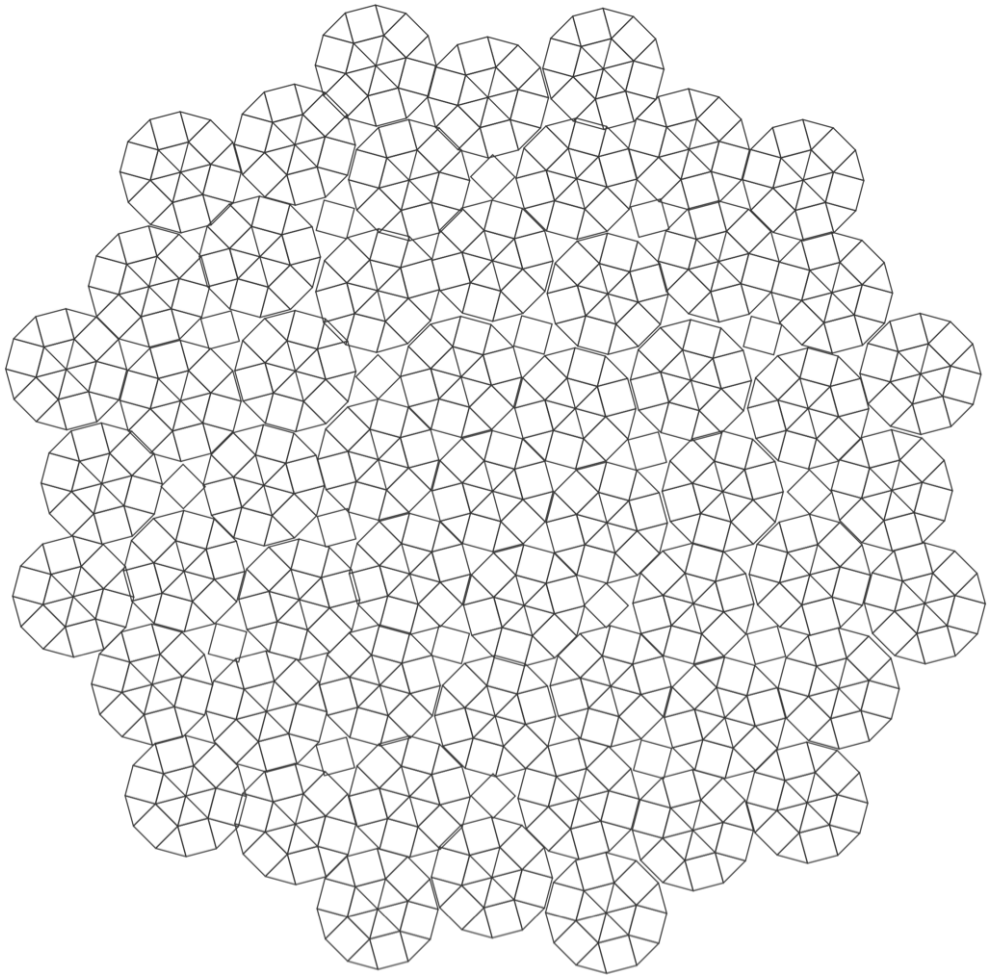

*Figure 3: The square-triangle tiling is made of two tiles: a square and an equilateral triangle.*

In (Van Dongen, 2021), we show that the square-triangle tiling is a special case of the rhombus-triangle tiling where the rhombus is a square. Here below are substitution rules for generating the square-triangle tiling. These rules are similar to the ones presented in (Van Dongen, 2021), although we simplified their pattern. Also, we present them here with no color as this is possible with such rules.

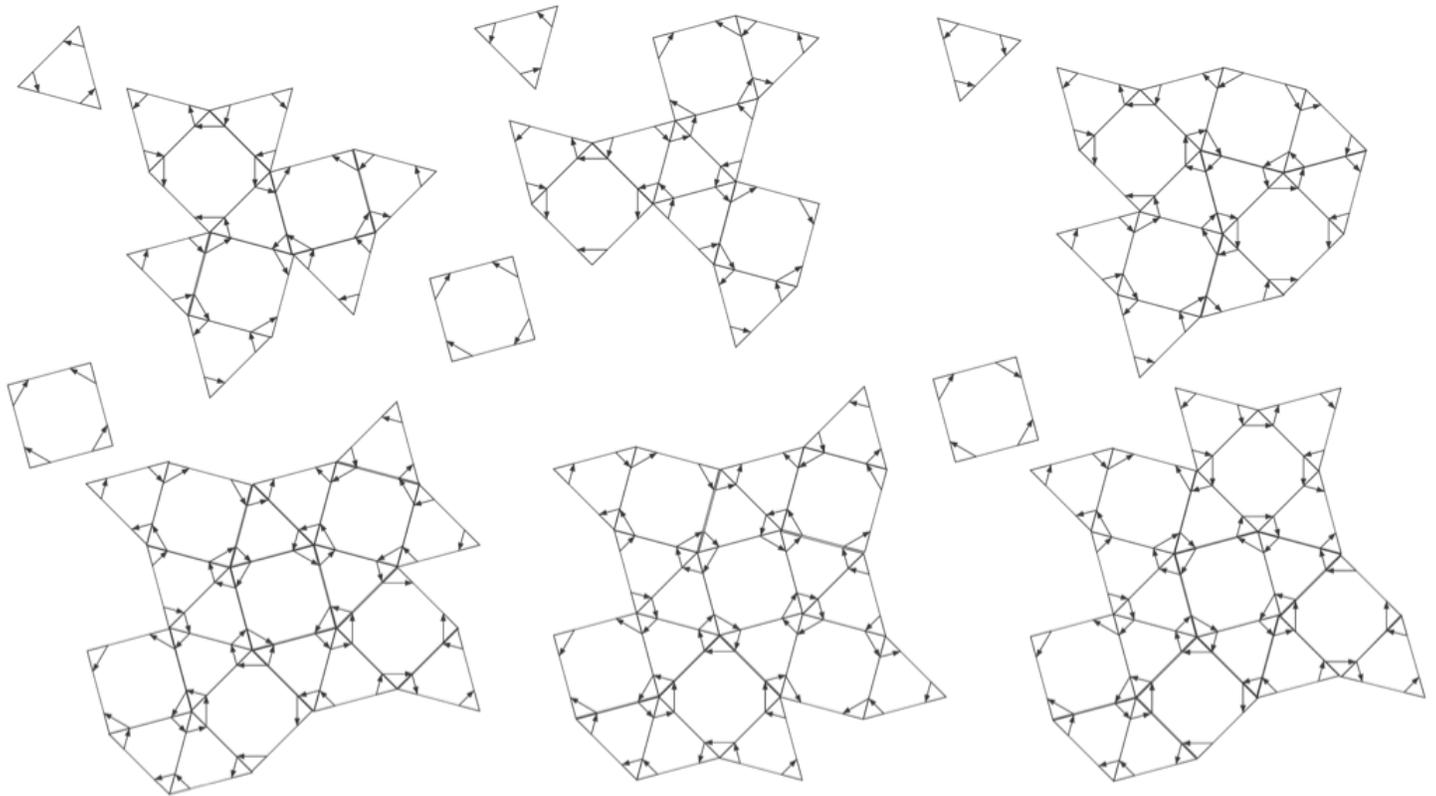

*Figure 4: The six rules of the rhombus-triangle tiling in the special case when the rhombus is a square.*

The matching rules on the tiles force the tiling to be aperiodic (and therefore structured). Here below is the application of these rules performed on the first tile. Note that substitution rules could be generated from any tile. Note also that the use of these rules allow for flipping and rotations of 60 degrees.

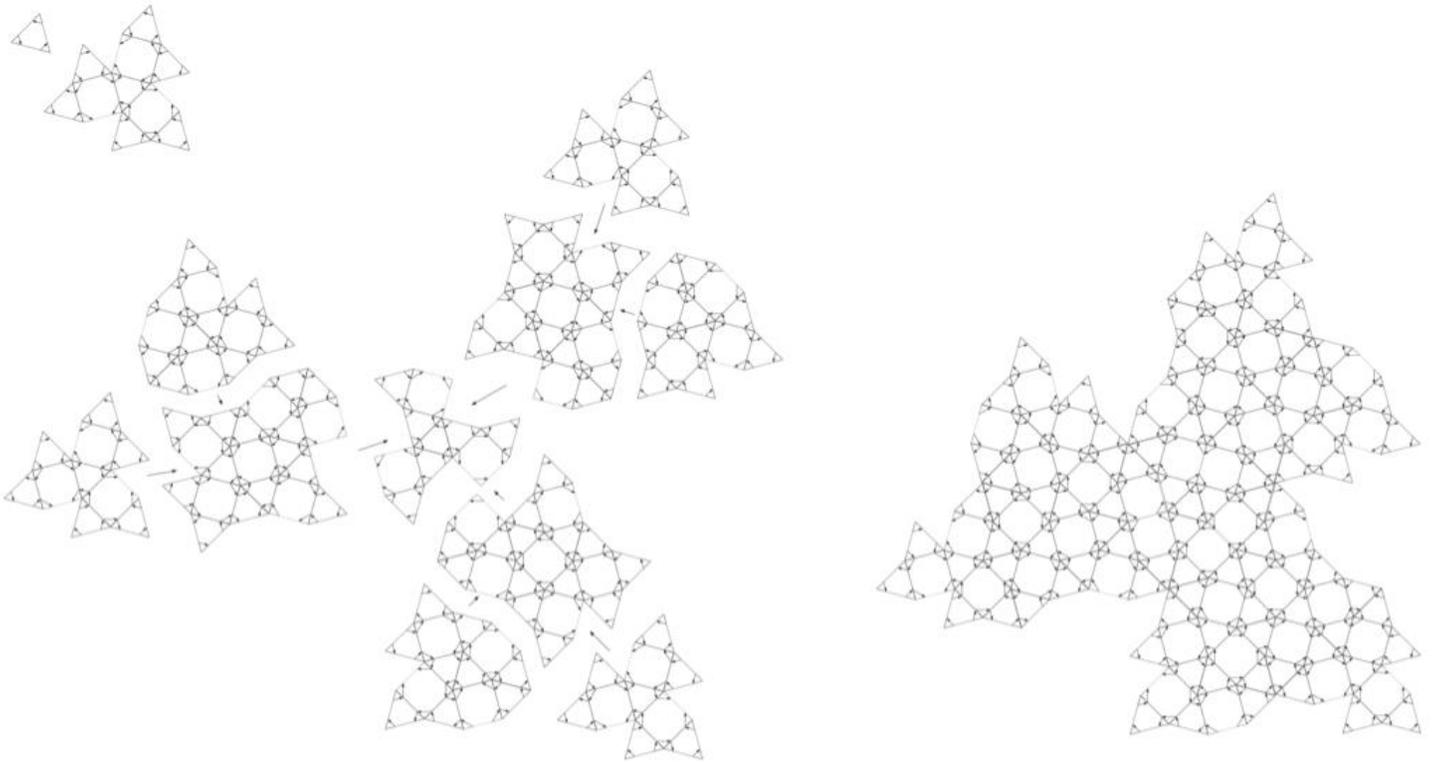

*Figure 5: Square-triangle aperiodic tiling: Example of use of the substitution rules applied to the first tile. In the middle, the tiles are being assembled. Please note that all arrows turn in one direction or another at any given vertex.*

Could there exist an aperiodic tiling made of one tile only? This is known as the ein-stein problem (Einstein Problem, s.d.). A first solution to this problem was proposed with a tile that is not connected. That is, it is not in a single piece. This proposal was made by Socolar-Taylor (Socolar, 2011). The tile is shown here below.

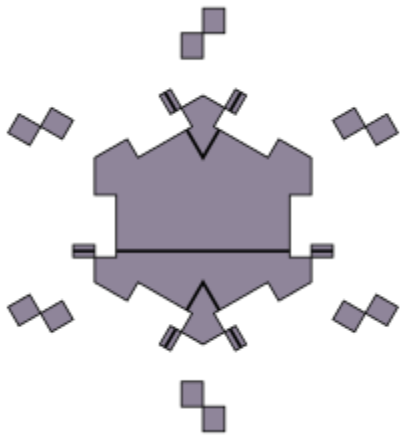

*Figure 6: Non-connected tile proposed as a solution to the ein-stein problem. Ref: https://en.wikipedia.org/wiki/Einstein_problem*

Another tiling was proposed recently (Walton & Whittaker, 2021). It is made of one hexagonal tile with a pattern on it. In this case, if one removes the pattern from the tile, aperiodicity disappears. The existence of a single connected tile, free of pattern and that can only cover the entire space in an aperiodic manner is still unknown until now. In this paper, we present an aperiodic tiling made of one tile: an isosceles right triangle. Patterns are needed to construct the aperiodic tiling but can be removed afterwards while keeping the tiling aperiodic.

## Substitution rules for the one-triangle aperiodic tiling

The aperiodic tiling made of one tile, an isosceles right triangle, can be constructed with fixed substitution rules; these are given here below.

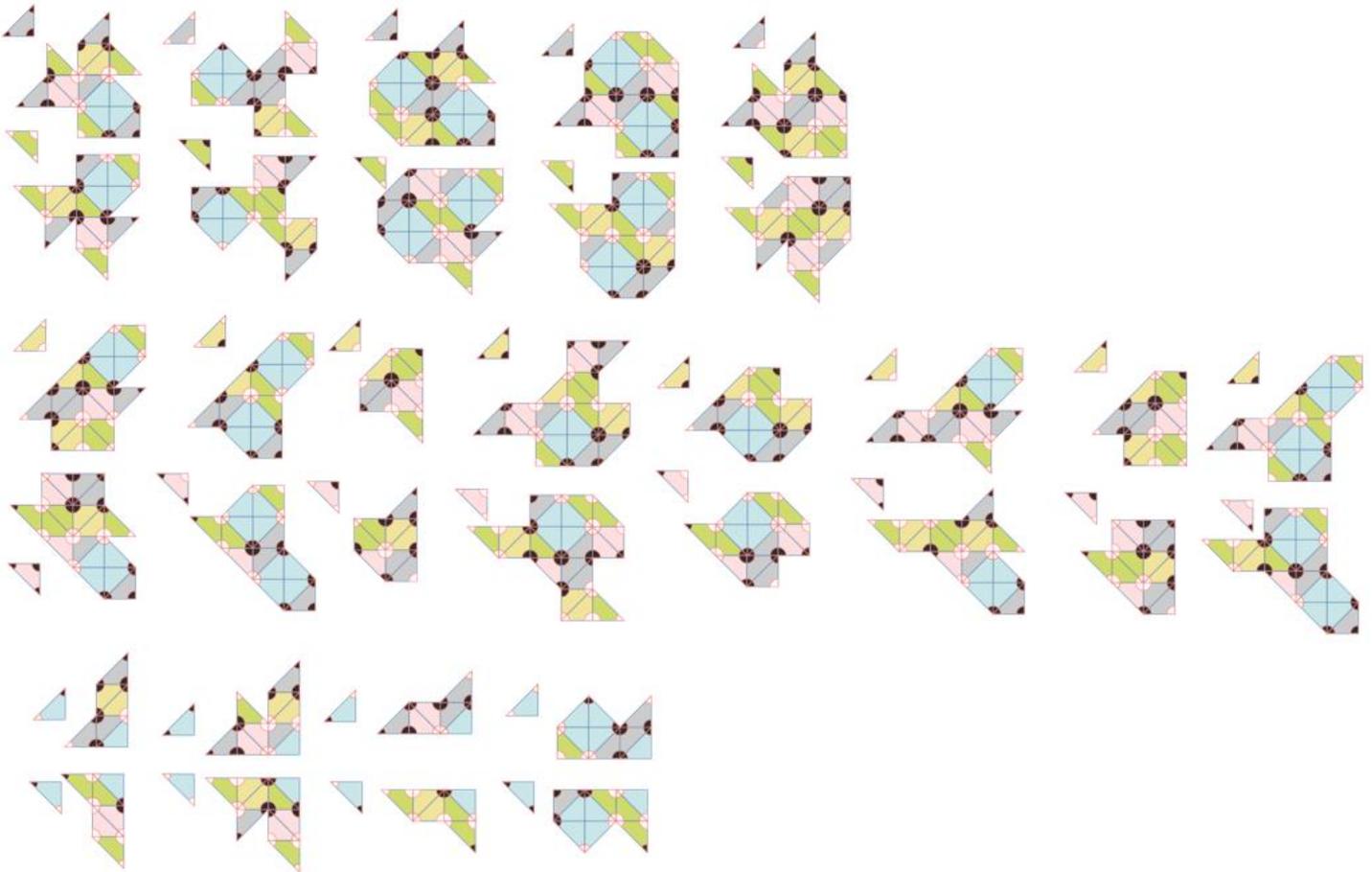

*Figure 7: Substitution rules to generate an aperiodic tiling made of one tile, a triangle.*

There are many rules compared to the square-triangle tiling but here the rotation of each tile is limited to 0 or 180 degrees, and no flipping is allowed. The rules are organized in pairs. There are 17 pairs of rules in total. Within each pair, one tile is the flip side of the other, with the following changes of colors:

- White corner switches to black corner and vice-versa;
- Green triangle tile switches to grey triangle tile and vice-versa;
- Yellow square tile switches to pink square tile and vice-versa.

Here is an example of application of these rules.

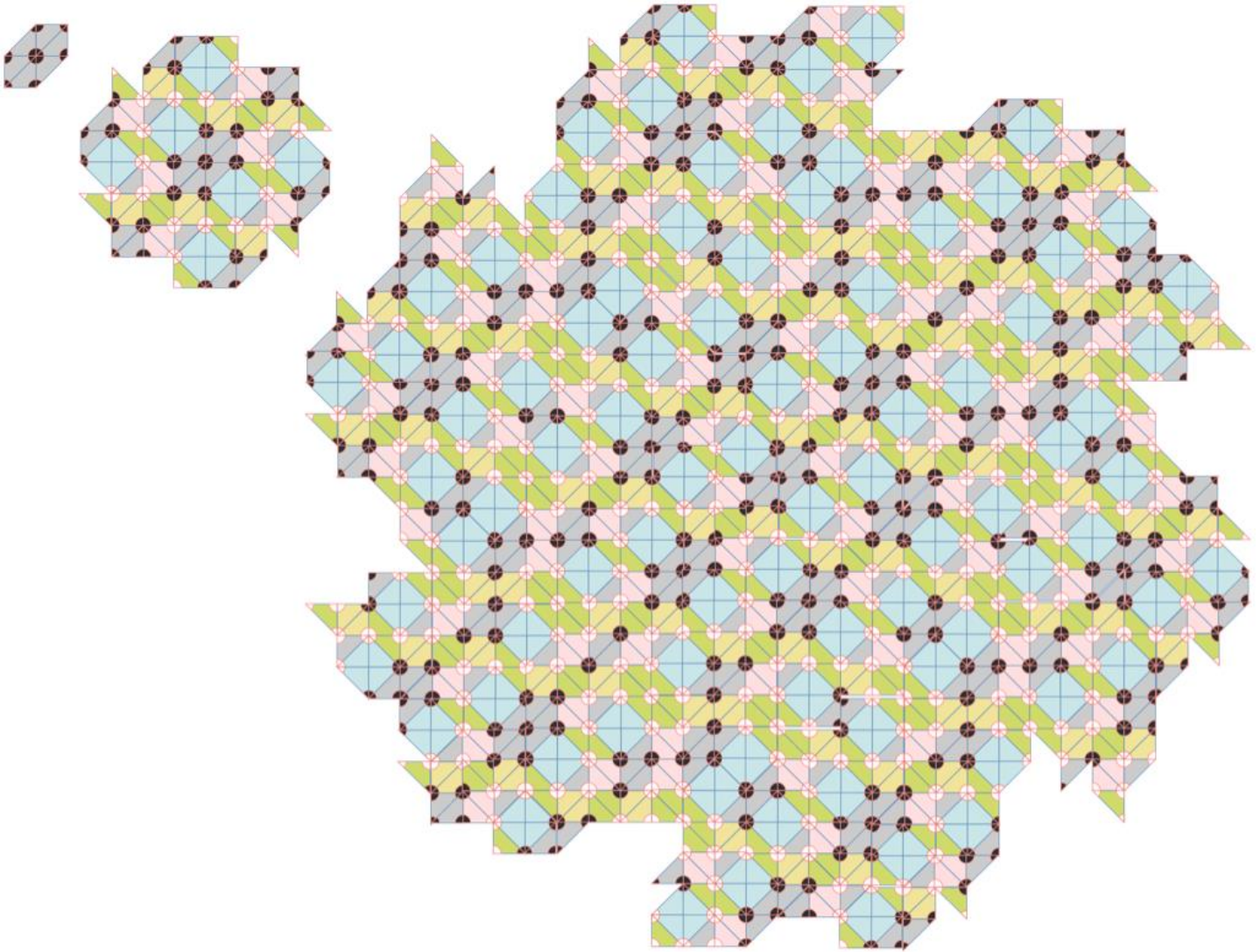

*Figure 8: Aperiodic one-triangle tiling: example of use of the substitution rules. All motifs can be removed once the tiling is created.*

# Proof of aperiodicity

There are several ways to prove that the proposed one-triangle tiling is aperiodic. We will prove it by comparing this new tiling to the square-triangle tiling that is known to be aperiodic. In the first part, we show that the one-triangle tiling is a direct reshaping of the square-triangle tiling. We present the conditions for the change of shape to work. In the second part, we show how all six substitution rules of the square-triangle tiling presented above can be rewritten into 28 substitution rules on two types of tiles: squares and isosceles right triangles. Finally, we show how to rewrite the rules on square tiles into rules on isosceles right triangles as well. This will complete the explanation of how to obtain the 17 pairs of rules presented here above. We will illustrate the case on some examples.

**A reshaping of tiles that works**

The following figure provides details on the mapping that we discovered to transform a square-triangle aperiodic tiling into our so-called one-triangle aperiodic tiling.

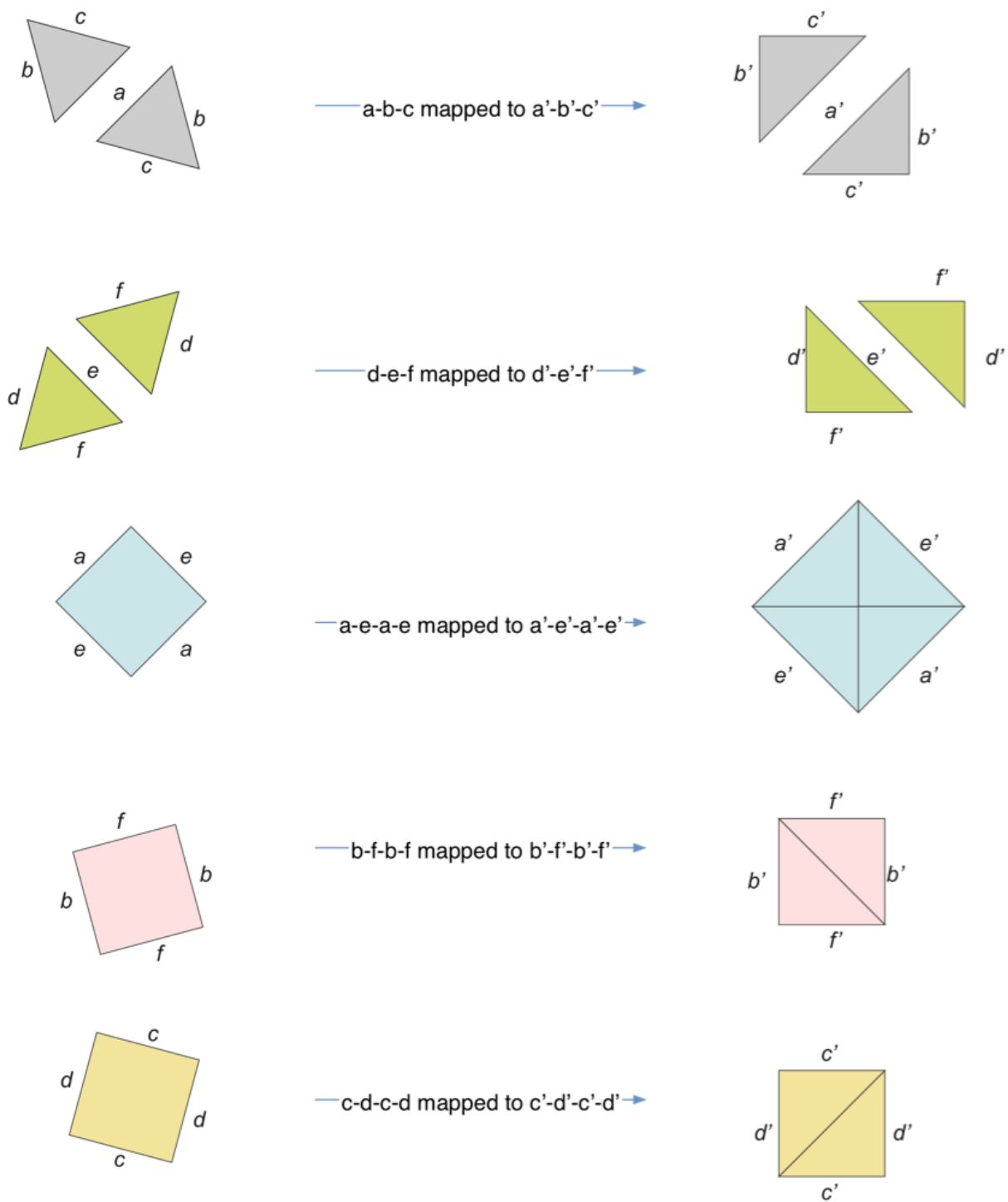

*Figure 9: On the left, all possible shapes and orientation of the tiles of the square-triangle tiling are shown (assuming that the tiling was properly oriented to exhibit some horizontal and vertical symmetry to start with). The mapping of each tile into new shapes is shown. On the right, all tiles consists of one or several isosceles right triangle, all of the same size. All edges that are of same length and same orientation are given a label (i.e. a, b, c, d, e, f). An edge of a particular label is mapped into a new edge that can be labeled the same way (i.e. a', b', c', d', e', f'). All edges labeled a are mapped to edges labeled a', etc. The transformation of the tiles follow the same logic. They are identified with the labels of their edges. Tile a-b-c is mapped into tile a'-b'-c', etc. On the left, tiles of the same label are given the same color; they are given the same color after their change of shape.*

On the left, all possible shapes and orientation of the tiles of the square-triangle tiling are shown (assuming that the tiling was properly oriented to exhibit some horizontal and vertical symmetry to start with). The mapping of each tile into new shapes is shown. On the right, all tiles consists of one or several isosceles right triangle, all of the same size. This was key in order for the new tiling to be of one tile. The transformation of the tiles is changing their shapes. This is done with strict rules so that edge-to-edge tiling is kept in place. The following rules are used:

- All edges that are of same length and same orientation are given a label (i.e. a, b, c, d, e, f). For example, all the 45-degree edges of unit length are given a label: <u>a</u>.
- An edge of a particular label is mapped into a new edge that can be labeled the same way (i.e. a', b', c', d', e', f'). All edges labeled <u>a</u> are mapped to edges labeled <u>a'</u>, etc. Assuming that <u>b'</u> is of unit length, the edge labeled <u>a'</u> is oriented at 45 degrees with length equal to sqrt(2).
- The transformation of the tiles follow the same logic. They are identified with the labels of their edges. Tile a-b-c is mapped into tile a'-b'-c', etc.
- On the left, tiles of the same label are given the same color; they are given the same color after their change of shape in order to facilitate the reading of the transformation.

The following figure illustrates the transformation. It starts with a square-triangle tiling with no color and ends with a one-triangle tiling with no color.

- Starting point: Square-triangle tiling with no color
- First transformation: Square-triangle tiling with the colors provided above
- Second transformation: One-triangle tiling with colors
- Final transformation: One-triangle tiling with all colors being removed

The following figure illustrates these steps.

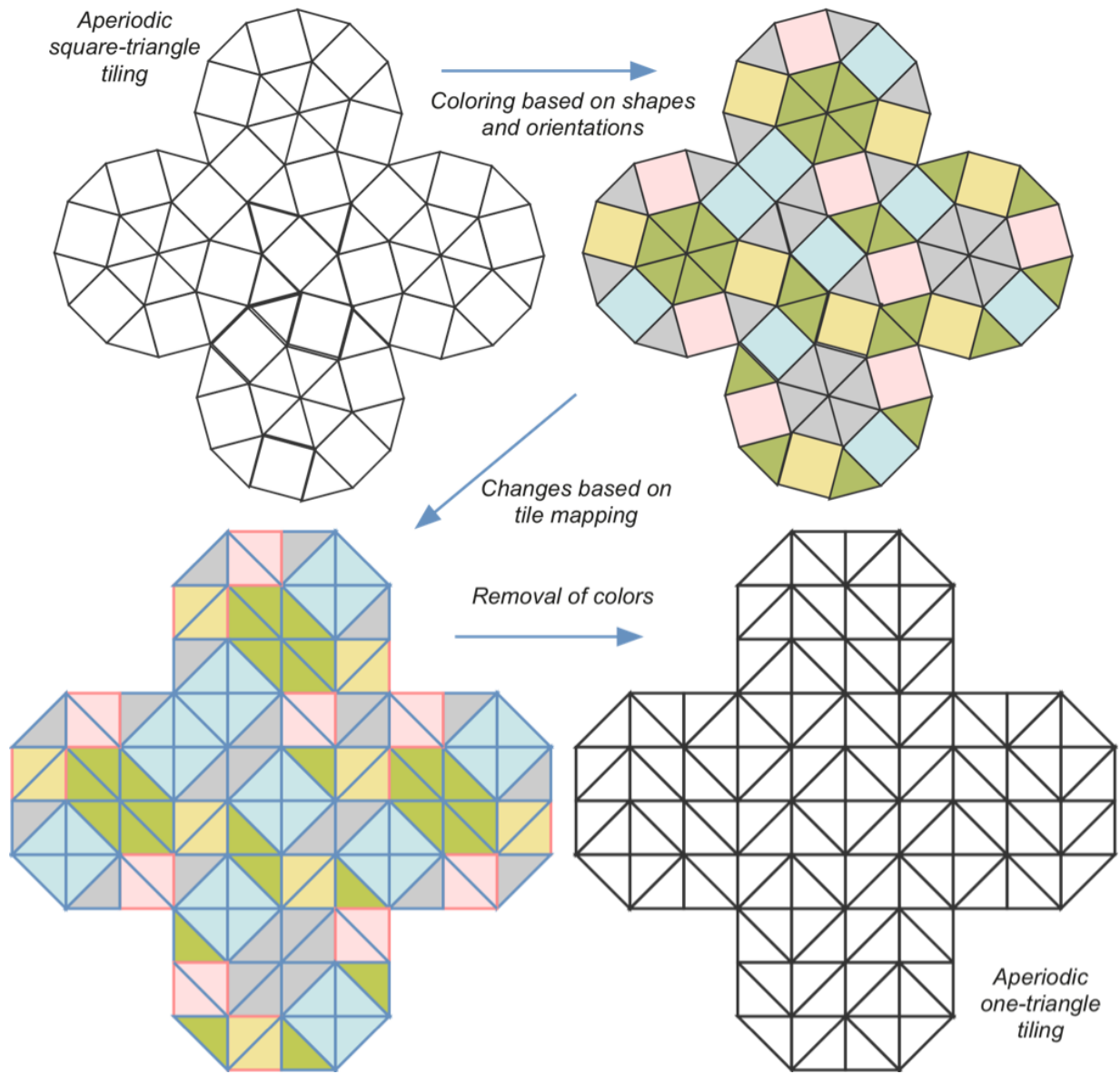

*Figure 10: Example of how a square-triangle aperiodic tiling is transformed into a one-triangle aperiodic tiling. It is done in three steps. First, the tiles of the square-triangle tiling are colored according to their shape and orientation. Second, the tiles are changed based on tile mapping. Third, the colors are removed to exhibit the fact that the new tiling is made of a single tile: an isosceles right triangle.*

**Development of the substitution rules for the one-triangle tiling**

The six substitution rules of the square-triangle tiling can be transformed using the method shown here above. However, no flipping and only rotations of either 0 or 180 degrees should be allowed on the tiles to keep the mapping valid. (Only a rotation of 180 degrees will preserve the labels given to the edges and shapes.)

We are now ready to start with the rewriting of the first substitution rule. The figure below shows how our first two new substitution rules for the one-triangle tiling are derived. Note the symmetry between the two new rules. The shapes are symmetrical but the tile colors differ.

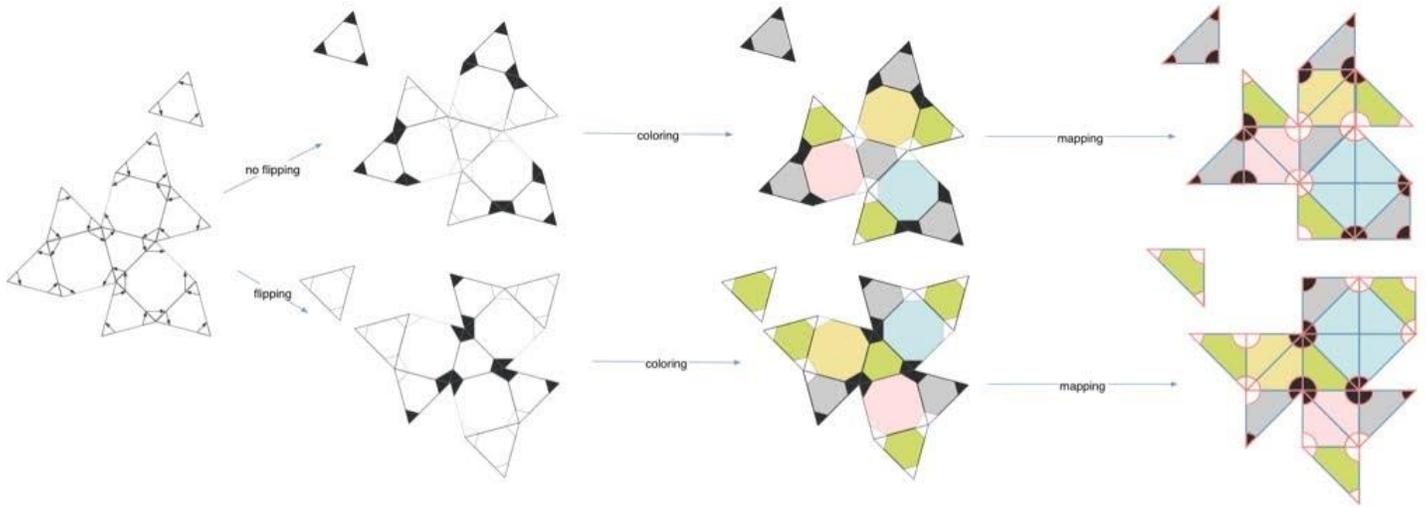

*Figure 11: Transformation of the first substitution rule of the square-triangle tiling into new rules for generating the one-triangle tiling. The first rule is mapped into two new rules. The use of these rules allow only rotations of 180 degrees, no flipping. The two new rules are symmetrical but with changes of colors.*

We can apply the same type of transformations on the next rule. The result is provided here below.

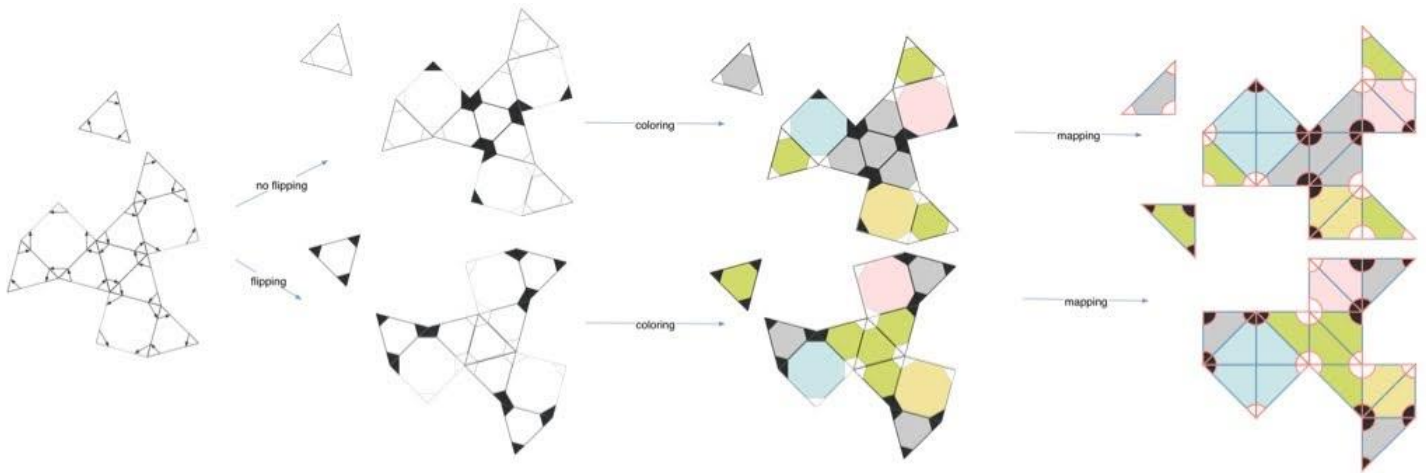

*Figure 12: Transformation of the second rule of the square-triangle tiling into two new rules for the one-triangle tiling.*

The first two rules were short to derive because rotations of 120 degrees and of 240 degrees keep them unchanged. The next four rules are different. They lead to the creation of 6 cases per rule instead of 2 as shown below for the next rule.

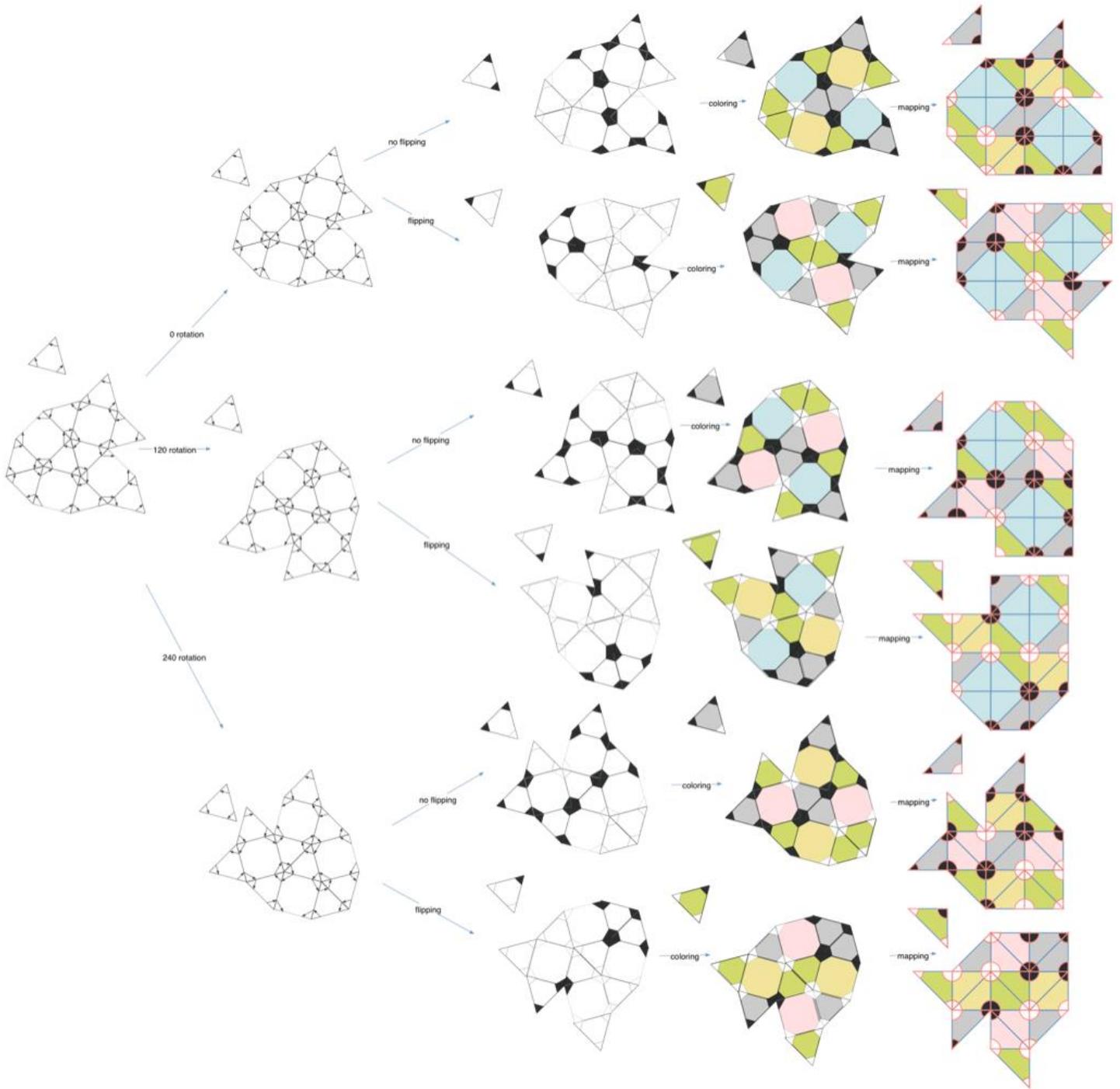

*Figure 13: Third rule transformed so that the resulting substitution rules only allow a rotation of either 0 or 180 degrees.*

Here below is the development of the fourth rule. The same method is being used in order to obtain substitution rules that only allow rotations of either 0 or 180 degrees and no flipping.

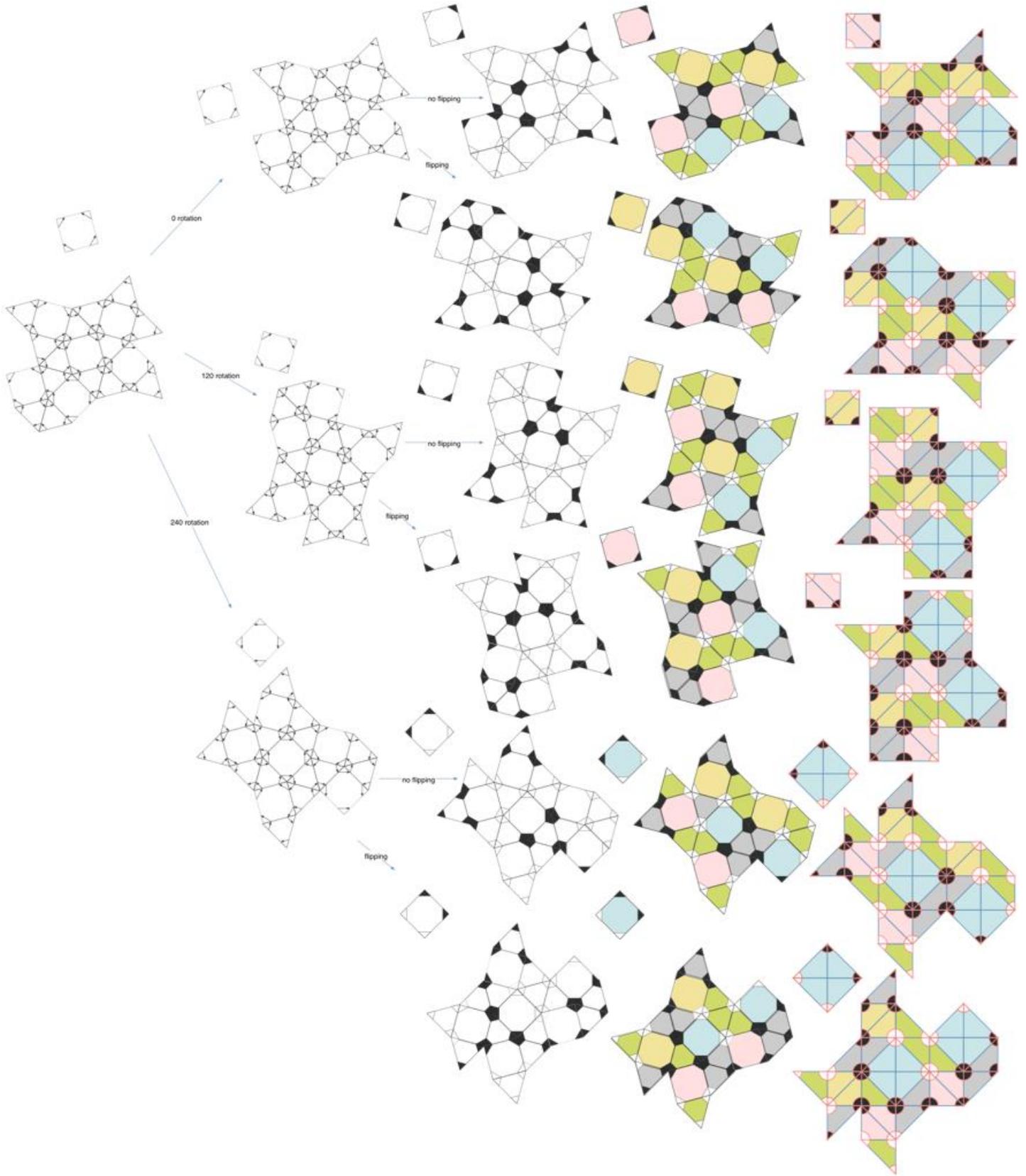

*Figure 14: Transformation steps on the fourth rule of the square-triangle tiling for the rules of the one-triangle tiling to allow only rotations of 0 or 180 degrees and no flipping.*

In the same manner, the fifth rule is transformed. The result is shown here below.

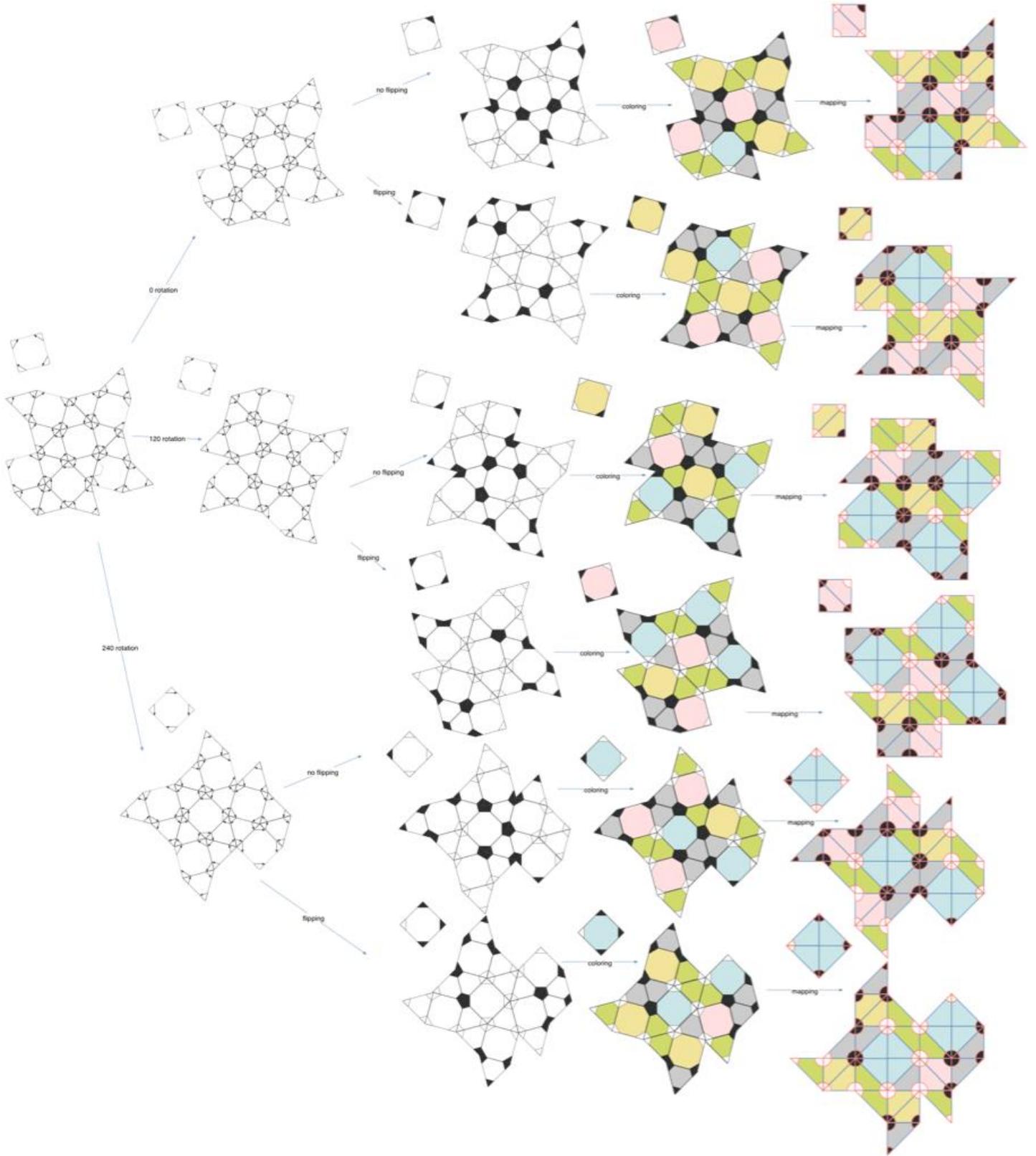

*Figure 15: Fifth rule of the square-triangle tiling transformed for the one-triangle tiling.*

Finally, the development of the last rule, using the same technique, is provided below.

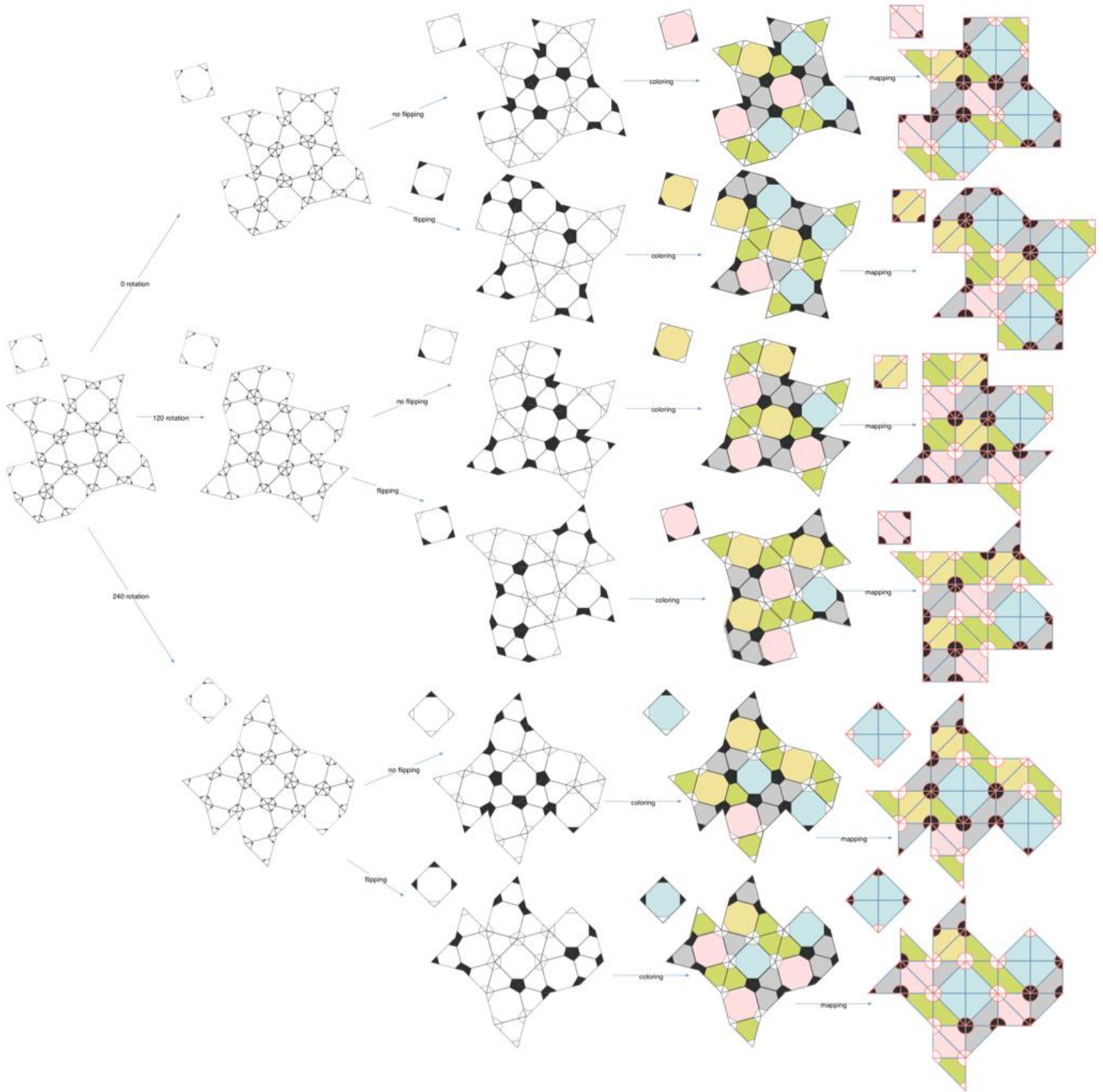

*Figure 16: Transformation of the sixth rule of the square-triangle tiling into substitution rules for the one-triangle aperiodic tiling.*

In order to complete this section, let us compare the use of the rules of the aperiodic square-triangle tiling and of the one-triangle one on the previous example.

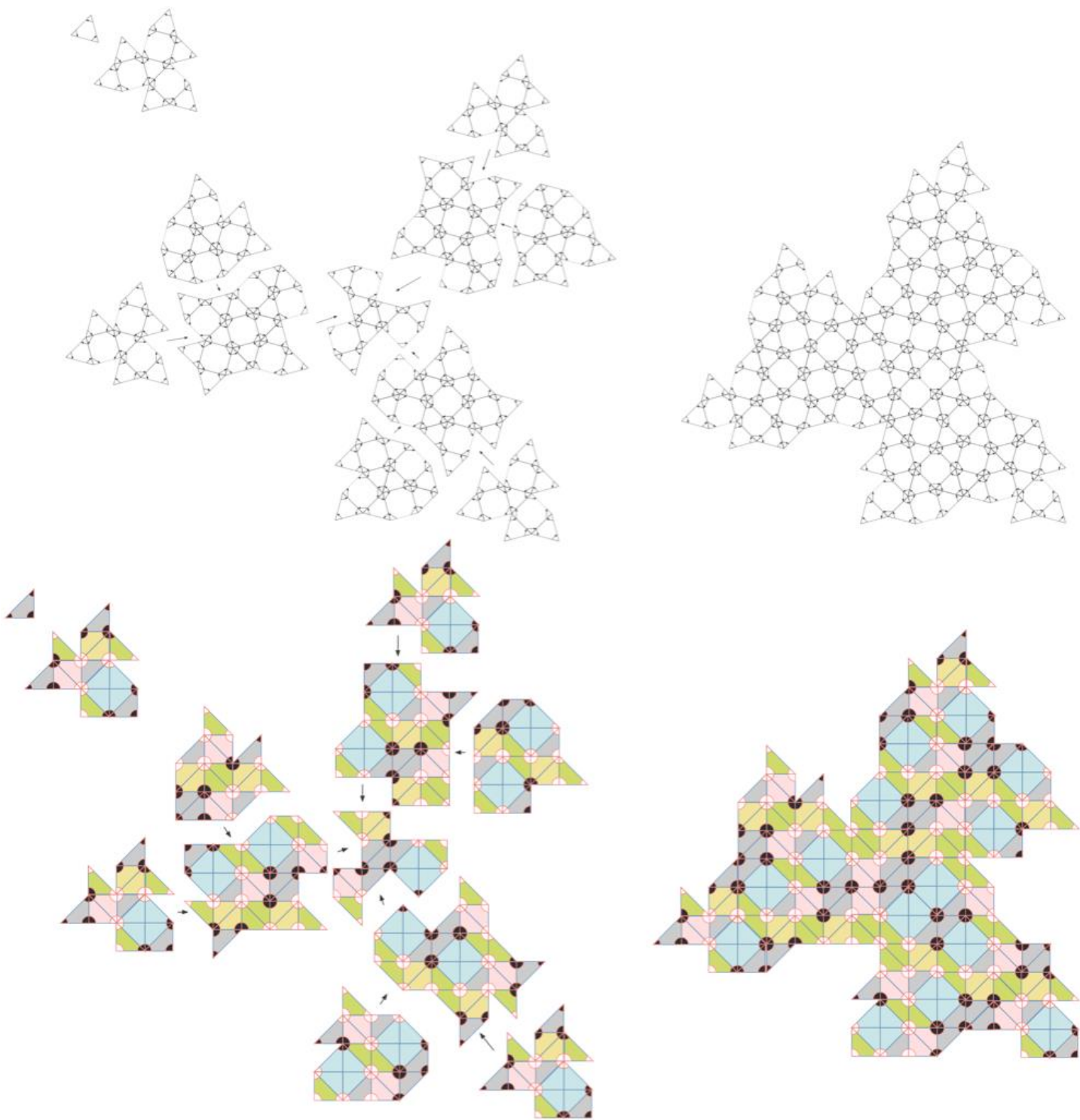

*Figure 17: Comparison on the use of the substitution rules between the aperiodic square-triangle tiling and the aperiodic one-triangle tiling.*

## Substitution rules on the same triangle only

The rules presented in the previous section make use of different tiles: isosceles right triangles and squares. Although the squares can be considered as a set of triangles, all substitution rules are not using the same tile. Therefore one could argue that the tiling is not made of one single tile. In this section, we show that this can be fixed. At the end, the substitution rules can all be on the same tile. This is done by cutting into pieces the substitution rules on the squares. Examples are shown here below. By applying the same technique on all the square-tile rules, all the rules can be using the same triangle: an isosceles right triangle. The new set of substitution rules is shown in Figure 7.

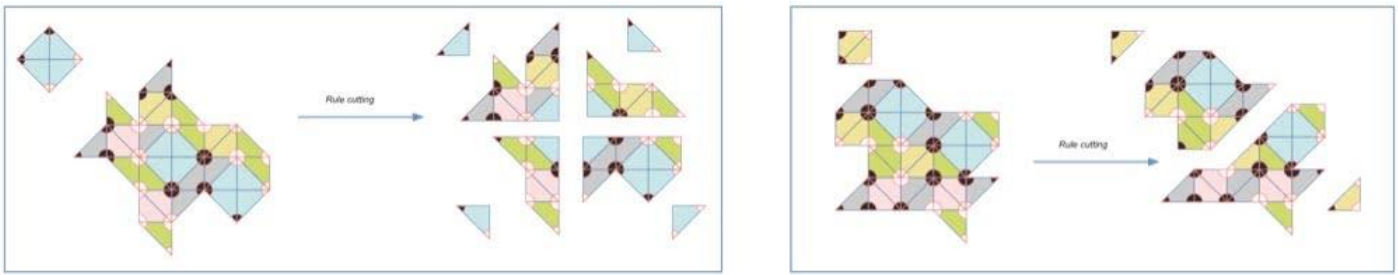

*Figure 18: Two examples of rule cutting: rules on squares are rewritten as rules on the same triangle. All motifs are kept the same.*

## Conclusion

In this paper, we presented a new aperiodic tiling made of a single tile, an isosceles right triangle. The tile itself is not aperiodic but the tiling is. It is created by means of substitution rules all on the same triangle but with different patterns. But once the tiling is generated, all patterns can be removed and the tiling remains aperiodic.

Because the one-triangle aperiodic tiling is based on the square-triangle one, it also possesses an underlying dodecagonal structure. Its use for the study of dodecagonal quasicrystals and soft matters in particular now needs to be investigated. The underlying structure of the aperiodic one-triangle tiling is dodecagonal but slightly irregular. Worth to be presented in details, it will be the subject of a separate document.

On top of presenting a new tiling, we presented in this paper a methodology for changing existing tilings and their associated substitution rules. In the future, it would be interesting to see whether the same method could be used to discover other new tilings.

## Special thanks

I would like to thank all the members of my family for being so supportive with my work. Also, I would like to thank Pierre Gradit, a researcher with whom I used to work some 30 years ago; his input was truly valuable.